\documentclass[12pt]{article}

\usepackage{color}
\usepackage{amsmath, amssymb}
\usepackage{amsthm,enumerate}

\newtheorem{assume}{Assumption}[subsection]

\newtheorem{theorem}{Theorem}[subsection]
\newtheorem{lemma}{Lemma} [subsection]
\newtheorem{prop}{Proposition} [subsection]

\newtheorem{rem}{Remark}[subsection]

\newcommand{\bec}[1]{\mbox{\boldmath $#1$}}

\allowdisplaybreaks

\begin{document}

\title{{\bf Limit theorems and bias correction for the initial-noise-adjusted SIML estimator}} 

\author{Rikuki Okamoto\footnote{
{\bf Corresponding author}: Rikuki Okamoto, Department of Mathematical Sciences, College of Science and Engineering, Ritsumeikan University, 1-1-1 Noji-higashi, Kusatsu, 525-8577, JAPAN \\
e-mail: ra0099vs@ed.ritsumei.ac.jp} \\
[8pt]
Ritsumeikan University
}
\date{\empty }

\maketitle

\begin{abstract}
We propose a sine-based modification of the separating information maximum likelihood (SIML) estimator for integrated volatility and covariance under equidistant sampling with microstructure noise. The estimator allows for nonzero initial observation noise. We establish consistency under moment assumptions and asymptotic normality under first-order Malliavin regularity of the diffusion coefficients. We also construct a feasible noise-bias correction using an estimator of the unknown noise covariance. This correction permits wider spectral bandwidths for consistency and stable asymptotic normality. 
\end{abstract}

\section{Introduction}

\subsection{High-frequency statistics under micro-structure noise: a literature review}

Let $(\Omega, \mathcal{F}, \mathbb{P},  \mathbb{F})$ be a complete filtered  probability space, which supports a $d$-dimensional Wiener process $\bec{W} = (W^1, W^2, ..., W^d)$. Here $\mathbb{F} :=  \{\mathcal{F}_t\}_{0 \leq t \leq 1}$ is the augmentation of the filtration generated by $\{\bec{W}_s\}_{0 \leq s \leq t}$. 
We consider a $ J $-dimensional It\^o process
 \begin{align*}
     X_t^j = X_0^j + \int_0^t b^j(s) ds + \sum_{r=1}^d \int_0^t \sigma^j_r(s) dW_s^r, 
     j=1,\cdots, J, 
 \end{align*}
where $b^j, \sigma^j_r \in L^2_a[0,1]$ for all $j=1,2,...,J$, $r=1,2,...,d$ and $t \in [0,1]$. Here we denote by $L^2_a[0,1]$ the space of $\mathbb{F} $-adapted process $\theta$ with 
the property $\mathbb{E}\left[ \int_0^1 \theta ^2(s) ds \right] < \infty$.

We suppose that we observe $ X^j $
at $\{t_k^j\}_{k=0}^{n^{j}} (0=t_0^j < t_1^j < \cdots < t_{n^j}^j = 1)$ for $ j=1, \cdots, J $.
In the literature of high-frequency financial econometrics, 
the focus is mostly on the estimation of the (integrated)
volatility 
\begin{align}\label{Volatility}
  \int_0^t \Sigma^{j,j'}(s) ds := \sum_{r=1}^d \int_0^t \sigma^j_r(s) \sigma^{j'}_r(s) ds, \quad  t \in (0,1],
\end{align}
relying on, 
instead of stationarity, the semi-martingale principle
\begin{equation*}
\sum_{k}(X_{t_{k}}-X_{t_{k-1}})^{2}\rightarrow[X]_{t} \quad \text{in probability as } \max_{k}|t_{k}-t_{k-1}|\rightarrow 0. 
\end{equation*}
However, in financial practice, market microstructure noise contaminates the efficient price; 
the latent process $X$ is not directly observable.
Instead, we observe noisy discrete-time data
\begin{equation}
  Y_{t_k} = X_{t_k} + v_{t_k}, \ k=0,1,\dots,n,
\label{eq:observation_model}
\end{equation}
where $(v_{t_k})$ represents microstructure noise.

If the noise term in \eqref{eq:observation_model} is not negligible, standard realized volatility and covariance estimators based on the increments
\begin{equation*}
  \Delta_i Y = Y_{t_i} - Y_{t_{i-1}}
\end{equation*}
become inconsistent (see e.g. \cite[Section 7.1.5]{AitSahaliaJacod2014}). 
This phenomenon reflects the fact that the quadratic variation of the noise dominates that of the latent process at high sampling frequencies.

A rigorous understanding of this problem was first provided by Gloter and Jacod \cite{GloterJacod2001a,GloterJacod2001b}.
In the parametric setting where $\Sigma_t \equiv \Sigma(\theta)$ (deterministic) depends on an unknown finite-dimensional parameter $\theta$,
they consider the joint distribution of the increment vector
\begin{equation*}
  \Delta Y = (\Delta_1 Y^\top,\dots,\Delta_n Y^\top)^\top, 
\end{equation*}
that is, under \eqref{eq:observation_model}
and 
\begin{equation}
 v_{t_i} \sim \mathcal{N}(0,\nu), 
 \qquad v_{t_i}\ \text{i.i.d.\ and independent of } X, 
 \label{eq:noise}
\end{equation}
it holds that 
\begin{equation}
  \Delta Y \sim \mathcal{N}\bigl(0,\, V_n(\theta)\bigr),
  \qquad 
  V_n(\theta)
  = \frac{1}{n}(I_n \otimes \Sigma(\theta))
    + (2I_n - J_+ - J_-) \otimes \nu,
  \label{eq:GJ_cov}
\end{equation}
where $I_n$ denotes the $n\times n$ identity matrix, 
and $J_+$ and $J_-$ as the first upper and lower shift matrices, respectively:
\begin{equation}
  (J_+)_{ij} = \mathbf{1}_{\{j=i+1\}},
  \qquad
  (J_-)_{ij} = \mathbf{1}_{\{j=i-1\}},
  \qquad i,j=1,\dots,n .
  \label{eq:Jpm_def}
\end{equation}
Based on this representation, Gloter and Jacod establish local asymptotic normality (LAN) of the statistical experiment and derive the asymptotic Fisher information.
A key implication of their theory is that, when the noise variance does not vanish with $n$,
the optimal rate of convergence for estimating $\theta$ is $n^{-1/4}$ rather than the classical $n^{-1/2}$ rate.
This slowdown is intrinsic to the noisy observation scheme and cannot be avoided by any estimator.

\medskip

The practical relevance of this LAN theory for financial econometrics was subsequently demonstrated by Xiu \cite{Xiu2010}.
Assuming additive Gaussian noise as in \eqref{eq:noise}, Xiu proposes a Gaussian quasi-maximum likelihood estimator (QMLE)
based on the criterion function
\begin{equation}
  \ell_n(\theta)
  = -\frac12 \Bigl(
      \log\det V_n(\theta)
      + \Delta Y^\top V_n(\theta)^{-1} \Delta Y
    \Bigr),
  \label{eq:QMLE}
\end{equation}
where $V_n(\theta)$ is given by \eqref{eq:GJ_cov}.
Even when the volatility process $\Sigma_t$ is stochastic and time-varying, Xiu shows that the resulting QMLE achieves
\begin{equation}\label{Xiuconvergence}
  n^{1/4}(\hat\theta_n - \theta)
  \;\Rightarrow\;
  \mathcal{MN}\bigl(0, I(\theta)^{-1}\bigr),
\end{equation}
with the same asymptotic information matrix $I(\theta)$ as in Gloter and Jacod.
This result establishes that likelihood-based methods can attain the efficiency bound predicted by the LAN theory in empirically relevant settings.

\medskip

An alternative likelihood-oriented approach is developed by Kunitomo, Sato and their collaborators through the separating information maximum likelihood (SIML) method \cite{KS08a, KS08b, KS10, KS11, KS13, KSK18}. 
Starting from the covariance structure \eqref{eq:GJ_cov}, SIML exploits the spectral decomposition of the matrix
\begin{equation*}
  2I_n - J_+ - J_-,
\end{equation*}
which diagonalizes the noise component of $V_n(\theta)$.
By transforming the data into the corresponding frequency domain, SIML separates low-frequency components, where the signal term $(1/n)\Sigma$ is dominant, from high-frequency components dominated by noise.
Rather than maximizing the full quasi-likelihood \eqref{eq:QMLE}, SIML constructs estimators using selected frequency components or weighted quadratic forms.

\medskip

In summary, Gloter and Jacod \cite{GloterJacod2001a,GloterJacod2001b} provide the fundamental asymptotic efficiency bounds for diffusion models observed with microstructure noise,
Xiu \cite{Xiu2010} demonstrates that these bounds are attainable by quasi-likelihood methods,
and the SIML framework of Kunitomo and Sato offers a practically flexible implementation based on spectral separation.

\subsection{The SIML estimator and  Malliavin-Mancino's Fourier method}
We continue to discuss SIML approach. 
A key issue in the SIML estimator construction is that they use 
\begin{align}\label{Jacobi1}
    J^{\text{SIML}} := 
    \begin{pmatrix}
        1 & 1 & 0 & \cdots & 0 & 0 & 0 \\
        1 & 0 & 1 & \cdots & 0 & 0 & 0 \\
        0 & 1 & 0 & \ddots & 0 & 0 & 0 \\
        \vdots & \vdots & \ddots & \ddots & \ddots & \vdots & \vdots \\
        0 & 0 & 0 & \ddots & 0 & 1 & 0 \\
        0 & 0 & 0 & \cdots & 1 & 0 & 1 \\
        0 & 0 & 0 & \cdots & 0 & 1 & 0
    \end{pmatrix}. 
\end{align}
instead of $ J_+ + J_- $ in defining the likelihood function. 
This comes from their assumption that $ v_{t_0} = 0 $. 
Noting that $ J^{\text{SIML}} $ is diagonalized by $ \mathbf{P}_n = (p^n_{k,l}) $ with
\begin{equation*}
    p^n_{k,l} = \sqrt{\frac{2}{n+ \frac{1}{2}}}
    \cos \left( \left(l - \frac{1}{2}\right)
    \pi \left( \frac{k-\frac{1}{2}}{n+\frac{1}{2}}
    \right) \right),
\end{equation*}
they work on 
the log likelihood function 
of $ \mathbf{z} = n^{1/2}
\mathbf{P}_n^{\top} {\Delta Y} $ assuming $ \Sigma(\theta) =c $, that is, 
\begin{align*}
    \begin{aligned}
        L_n (c, \nu) 
       := -\frac{1}{2} \sum_{k=1}^n \log (a_{k,n} \nu + c) -\frac{1}{2} \sum_{k=1}^m
       \frac{\mathbf{z}_k^2}{a_{k,n}\nu + c},    
    \end{aligned}
\end{align*}
where we put
\begin{align*}
    \begin{aligned}
        a_{k,n} := 4 n \sin^2 \left(\frac{2k-1}{2n+1}\frac{\pi}{2}  \right), \quad k=1, 2, \dots, n. 
    \end{aligned}
\end{align*}

Another key issue in SIML method is the decomposition 
\begin{align*}
    \begin{aligned}
        2 L_n = L_n^{(1)} + L_n^{(2)} + L_n^{(r)},
    \end{aligned}
\end{align*}
where we put
\begin{align}\label{KSL0}
    \begin{aligned}
        L_n^{(1)} (c) 
        := -m \log c - \frac{1}{c}\sum_{k=1}^m \mathbf{z}_k^2,  
    \end{aligned}
\end{align}
\begin{align*}
    \begin{aligned}
        L_n^{(2)} (\nu) 
        := - \sum_{k=n+1-l}^n \log (a_{k,n} \nu ) -\frac{1}{\nu} \sum_{k=n+1-l}^n \frac{\mathbf{z}_k^2}{a_{k,n}},
    \end{aligned}
\end{align*}
and $ L_n^{(r)} := 2 L_n - L_n^{(1)} - L_n^{(2)}$.
Note that the decomposition depends on the choice of $ (m,l) $. 
The term {\it separating information} comes from this decomposition. 

The SIML estimator is given by the maximizer 
of $ L_n^{(1)} $, that is,  
\begin{equation}\label{SIMLestimator}
    \begin{split}
        V_{n,m_n}^{j,j'}
        := \frac{n}{m_n} \sum_{l=1}^{m_n} \left( \sum_{k=1}^{n^j}
        p^{n^j}_{k,l} \Delta Y^j_k \right)
        \left(\sum_{k'=1}^{n^{j'}} p^{n^{j'}}_{k',l} \Delta Y^{j'}_{k'} \right),  \qquad j, j'=1, 2, \dots, J. 
    \end{split}
\end{equation}

They have proved the following two properties:
\begin{enumerate}[{\bf (a)}]
    \item ({\it the consistency}):
    the convergence in probability of $ V^{j,j'}_{n, m_n} $ to $ \int_0^1 \Sigma^{j,j'} (s) \,\mathrm{d}s $ as $ n \to \infty $
    is attained, provided 
    that $ m_n = o (n^{1/2})$, and
    \item ({\it the asymptotic normality of the error}): the stable convergence of 
    \begin{equation*}
    \begin{split}
       & \sqrt{m_n}\left(V^{j,j'}_{n, m_n}-\int_0^1 \Sigma^{j,j'} (s) \,\mathrm{d}s \right) \\
        & \to N \left(0, \int_0^1 \left(\Sigma^{j,j}(s) \Sigma^{j',j'}(s) + (\Sigma^{j,j'} (s))^2 \right) \mathrm{d}s
        \right)
        \end{split}
    \end{equation*}
    holds true as $ n \to \infty$ 
    if $ m_n = o (n^{2/5})$. 
\end{enumerate}
under the following conditions (see \cite{KSK18} for details):
\begin{enumerate}[(i)]
    \item the observations 
are equally spaced, that is, $ t^j_k \equiv k/n $, 
\item $ \{v^j_k\}_{k=1}^n $ are
independently and identically distributed, and the common distribution is independent of the choices of $ n $ and $ j $,
\item $ v_0^j \equiv 0 $, 
\item $ \sigma $ is deterministic and $ b \equiv 0 $. 
\end{enumerate}

\begin{rem}
In \cite{ANW2023}, the consistency and the asymptotic normality are proven when $ b $ and $ \Sigma $ are adapted 
processes with some regularity under a general sampling scheme, but without micro-structure noise. 
\end{rem}

\begin{rem}
An estimator of $ \nu $ is also given by the maximizer of $ L_n^{(2)} $: 
\begin{align*}
    \nu^*_{n,l} := \frac{1}{l} \sum_{k=n+1-l}^n \frac{\mathbf{z}_k^2}{a_{k,n}}.
\end{align*}
\end{rem}

In \cite{ANW2024}, the assumption that $ v_0 = 0 $ is found to be critical; 
the SIML estimator \eqref{SIMLestimator} fails to be a consistent estimator
of the integrated volatility in the presence of  $ v_0 $. 
Besides this observation,  
\cite{ANW2024} reports that 
Malliavin-Mancino's Fourier based estimator \cite{MMFS, MMAS}
of integrated volatility is characterized as the maximizer of $ L^{(1)} $ under 
\begin{align}\label{JacobiMM}
    {J}^{\text{MM}} := 
    \begin{pmatrix}
        0 & 1 & 0 & \cdots & 0 & 0 & 1 \\
        1 & 0 & 1 & \cdots & 0 & 0 & 0 \\
        0 & 1 & 0 & \ddots & 0 & 0 & 0 \\
        \vdots & \vdots & \ddots & \ddots & \ddots & \vdots & \vdots \\
        0 & 0 & 0 & \ddots & 0 & 1 & 0 \\
        0 & 0 & 0 & \cdots & 1 & 0 & 1 \\
        1 & 0 & 0 & \cdots & 0 & 1 & 0
    \end{pmatrix}, 
\end{align}
instead of $ J^{\text{SIML}} $, and observes that 
it is a consistent estimator only when $ v_0 = v_{t_n} $. 

Given these observations, \cite{ANW2024} proposes an estimator that maximizes 
$ L^{(1)} $ under $ J_+ + J_- $ in place of $ J^{\text{SIML}} $, and conjectures that it is a consistent estimator under the standard assumption \eqref{eq:noise}, 
The present paper proves the conjecture and establishes the asymptotic normality 
under common equidistant sampling. 

\subsection{Alternative estimation methods to the SIML approach}
As a method to overcome the difficulties described above, a modification of the SIML estimation approach is suggested in \cite{ANW2024}. Specifically, the matrix in the original SIML method $J^{\text{SIML}}$ is replaced with the following matrix: 
\begin{align}
    \tilde{J}^{\text{SIML}} := 
    \begin{pmatrix}
        0 & 1 & 0 & \cdots & 0 & 0 & 0 \\
        1 & 0 & 1 & \cdots & 0 & 0 & 0 \\
        0 & 1 & 0 & \ddots & 0 & 0 & 0 \\
        \vdots & \vdots & \ddots & \ddots & \ddots & \vdots & \vdots \\
        0 & 0 & 0 & \ddots & 0 & 1 & 0 \\
        0 & 0 & 0 & \cdots & 1 & 0 & 1 \\
        0 & 0 & 0 & \cdots & 0 & 1 & 0
    \end{pmatrix}.
\end{align}
Here, the matrix $\tilde{J}^{\text{SIML}}$ introduced in \cite{ANW2024} can be expressed as $\tilde{J}^{\text{SIML}}=J_+ + J_-$.
In \cite{ANW2024}, the diagonalization of  $\tilde{J}^{\text{SIML}}$ is shown to be established in the following form :
\begin{lemma}[\cite{ANW2024}]
    Let $n \times n$ matrix $R_n$ as 
    \begin{align*}
        r_{\ell,k}^{n} := \sqrt{\frac{2}{n+1}}\sin{\left( \frac{k\ell \pi}{n+1}\right)}, \ell,k= 1,2,..., n.
    \end{align*}
 
    Then, $R_n$ is the orthogonal matrix and 
    \begin{align}
        \tilde{J}^{SIML} = R_n \text{diag}\left[ 2\cos\left( \frac{\pi}{n+1} \right), 2\cos\left( \frac{2\pi}{n+1} \right), ..., 2\cos\left( \frac{n\pi}{n+1} \right) \right]  R_n^\mathsf{T} .
    \end{align}
\end{lemma}
We construct an estimation function from the components of this diagonalized matrix to obtain the integrated volatility that is robust against initial noise as follows; 
\begin{align}\label{modified_estimator}
    \hat{V}_{n,m_n}^{j,j'} &:= 
    \frac{n+1}{m_n}
    \sum_{l=1}^{m_n}
    \left( \sum_{k=1}^{n^j}r_{k,l}^{n^j}\Delta Y_k^j\right)
    \left( \sum_{{k'}=1}^{n^{j'}}r_{{k'},l}^{n^{j'}}\Delta Y_{k'}^{j'}\right). 
\end{align}

We refer to this estimator as the
\textit{Initial-noise-adjusted SIML estimator}.

In this paper, we study the initial-noise-adjusted SIML estimator
under common equidistant sampling with additive microstructure noise.
The noise vectors are independent across observation times and
independent of the signal, while dependence across components is
allowed. Under the stated coefficient and noise moment assumptions,
we establish consistency for integrated volatility and covariance when
$m_n \to \infty$ and $m_n=o(n^{1/2})$.
Under additional first-order Malliavin regularity of the diffusion
coefficients, we prove a centered stable mixed-normal limit at the
$\sqrt{m_n}$ scale when $m_n=o(n^{2/5})$.
The bracket convergence, asymptotic orthogonality, and negligibility
of the signal remainders are verified in the proof.

We also construct a feasible bias correction for unknown noise
covariance. The uncorrected pure-noise mean is of order $m_n^2/n$
when the corresponding noise covariance is nonzero, which explains
the restriction $m_n=o(n^{2/5})$ for the centered limit.
We estimate the noise covariance with error $O_p(n^{-1/2})$
and use it to correct this mean.
For both the oracle and feasible corrected estimators, consistency
holds when $m_n\to\infty$ and $m_n=o(n^{2/3})$.
Under the same additional Malliavin assumption, a centered stable
mixed-normal limit with the signal variance holds when
$m_n\to\infty$ and $m_n=o(n^{1/2})$.

\subsection{Organization of the present paper}

The remainder of the paper is organized as follows.
Section $2$ specifies the common equidistant sampling scheme and
the coefficient and noise assumptions.
It establishes the sine-difference and orthogonality identities,
decomposes the estimator into signal, cross, and pure-noise terms,
and proves consistency of the uncorrected estimator.

Section $3$ proves stable asymptotic normality under first-order
Malliavin regularity.
Deterministic kernel estimates, concentration near the diagonal,
and Malliavin duality are used to control the signal remainders
and verify bracket convergence and asymptotic orthogonality.
The signal limit is then extended to the uncorrected noisy
estimator.

Section $4$ constructs oracle and feasible noise-bias corrections,
establishes the accuracy of the noise covariance estimator, and
proves consistency and stable asymptotic normality of the corrected
estimators.
Section $5$ concludes the paper.

\section{Consistency}\label{sec:consistency}
\subsection{Sampling convention and assumptions}
Throughout Sections 2--5 we use the common equidistant scheme $n^j=n$ and $t_k^j \equiv t_k =k/(n+1)$, $k=0,\ldots,n$. For notational convenience, set $t_{n+1}=1$. This endpoint is auxiliary: the estimator uses only the $n$ observed increments indexed by $k=1,\ldots,n$, and no additional observation at $t_{n+1}$ is introduced. In particular, $t_n^j<1$; the unobserved terminal interval has length $1/(n+1)$.  This specifies the indexing consistently without changing the estimator.

Write $\Delta Y_k^j := Y_{t_k}^j-Y_{t_{k-1}}^j$, $\Delta X_k^j := X_{t_k}^j-X_{t_{k-1}}^j$, and $\Delta(v_k^j) := v_k^j-v_{k-1}^j$. The $J$ components are driven by the $d$-dimensional Wiener process from Section 1. Thus component indices run from $1$ to $J$, and Wiener indices from $1$ to $d$.

For every finite $p\ge2$, assume that
\begin{equation*}
A_p:=\mathbb E\left[\left(\sup_{0\le t\le1}\sum_j|b^j(t)|^2\right)^{p/2}\right]
+\mathbb E\left[\left(\sup_{0\le t\le1}\sum_{j,r}|\sigma_r^j(t)|^2\right)^{p/2}\right]<\infty. 
\end{equation*}
The processes are progressively measurable, and $\sigma$ has continuous sample paths. Let $(v_k^1,\ldots,v_k^J)$, $k\ge0$, be centered i.i.d. vectors independent of $\mathcal F_1$, with finite fourth moments. Write
$\nu^{j,j'}=\mathbb E[v_k^jv_k^{j'}]$.
The second-moment estimates in Sections 2 and 4 require only $A_4<\infty$.
The higher-order moment assumptions are used in the Malliavin estimates of Section 3. We take $\mathcal F_t$ to be the usual augmentation of $\sigma(W_s:0\le s\le t)$, and assume that the signal increments are measurable with respect to $\mathcal F_1$. 

To rewrite the finite-sum representation of the estimator in terms of stochastic integrals, we define the following function $\phi_\ell$:
\begin{equation*}
    \phi_\ell(u)=\begin{cases}
    \frac{\pi\ell k} {n+1} , & u\in(t_{k-1}^j,t_k^j],\ 1\le k\le n,\\ 0, &u\in(t_n^j,1] \text{ or }u=0.\end{cases}
\end{equation*}

At observation points the sine factors remain exactly those in Section 1. Define
\begin{equation*}
\mathcal D_{n,m_n}^{j,j'}(u,s)=\frac2{m_n}\sum_{\ell=1}^{m_n}\sin(\phi_\ell(u))\sin(\phi_\ell(s)).
\end{equation*}
We denote the latent-process (signal) estimator by $S_{n,m_n}^{j,j'}$:
\begin{align*}
S_{n,m_n}^{j,j'}
&=\frac{n+1}{m_n}\sum_{\ell=1}^{m_n}
\left(\sum_{k=1}^nr_{k,\ell}^n\Delta X_k^j\right)
\left(\sum_{k'=1}^nr_{k',\ell}^n\Delta X_{k'}^{j'}\right)\\
&=\frac2{m_n}\sum_{\ell=1}^{m_n}
\left(\int_0^1\sin(\phi_\ell(u))dX_u^j\right)
\left(\int_0^1\sin(\phi_\ell(s))dX_s^{j'}\right).
\end{align*}
The exact decomposition of the Initial-noise-adjusted SIML estimator is
\begin{equation*}
\hat V_{n,m_n}^{j,j'}=S_{n,m_n}^{j,j'}+C_{n,m_n}^{j,j'}+N_{n,m_n}^{j,j'},
\end{equation*}
where
\begin{align*}
C_{n,m_n}^{j,j'}&=\int_0^1\sum_{k'=1}^n\mathcal D_{n,m_n}^{j,j'}(s,t_{k'}^{j'})\Delta(v_{k'}^{j'})dX_s^j
+\int_0^1\sum_{k=1}^n\mathcal D_{n,m_n}^{j,j'}(t_k^j,s)\Delta(v_k^j)dX_s^{j'},\\
N_{n,m_n}^{j,j'}&=\sum_{k,k'=1}^n\mathcal D_{n,m_n}^{j,j'}(t_k^j,t_{k'}^{j'})\Delta(v_k^j)\Delta(v_{k'}^{j'}).
\end{align*}
To interpret the stochastic integrals involving the observation noise, we use the initially enlarged filtration $\mathcal F_t^v:=\mathcal F_t\vee\sigma(v_k:k\ge0)$, with the usual augmentation. Independence of the noise from $\mathcal F_1$ ensures that \(W\) remains a Brownian motion in this filtration. The cross-term stochastic integrals are understood relative to $\mathbb{F}^v = \{\mathcal F_t^v\}_{0 \leq t \leq 1}$, whereas the signal martingale and Malliavin arguments are carried out in the original filtration $\mathbb{F}$.

\subsection{Discrete orthogonality and the difference lemma}

\begin{lemma}\label{keylemma}
For $1\le\ell,\ell'\le m_n\le n$, the signed difference sum is
\begin{align*}
&\sum_{k=1}^n\Delta_{k}\left(\sin\frac{\pi\ell k}{n+1}\right)
\Delta_{k}\left(\sin\frac{\pi\ell' k}{n+1}\right)\\
&=2(n+1)\sin^2\frac{\pi\ell}{2(n+1)}\mathbf1_{\{\ell=\ell'\}}
-(-1)^{\ell+\ell'}\sin\frac{\pi\ell}{n+1}\sin\frac{\pi\ell'}{n+1}.
\end{align*}
Consequently its absolute value is $O(\ell^2/n)$ on the diagonal and $O(\ell\ell'/n^2)$ off the diagonal, with constants uniform in the indices. If $m_n=o(n)$, the diagonal is asymptotic to $\pi^2\ell^2/(2(n+1))$, uniformly in $\ell\le m_n$.
\end{lemma}

\begin{proof}
Fix $\ell,\ell'\in\{1,2,\ldots,m_n\}$. The sum-to-product identity gives
\begin{equation*}
\Delta_{k}\left(\sin\frac{\pi\ell k}{n+1}\right)
=2\sin\frac{\pi\ell}{2(n+1)}
\cos\frac{\pi\ell(k-\tfrac12)}{n+1}.
\end{equation*}
We first verify the finite cosine sums used below. For a nonzero integer $q$ with $ 0 < |q|<2(n+1)$, the identity
\begin{equation*}
2\sin\frac{\pi q}{2(n+1)}
\cos\frac{\pi q(k-\tfrac12)}{n+1} =\sin\frac{\pi q k}{n+1}
-\sin\frac{\pi q(k-1)}{n+1}
\end{equation*}
can be summed over $k=1,\ldots,n+1$. Since the right-hand side telescopes to $\sin(\pi q)=0$, and the factor on the left is nonzero, we obtain
\begin{equation*}
\sum_{k=1}^{n+1}\cos\frac{\pi q(k-\tfrac12)}{n+1}=0, \quad (q \neq 0).
\end{equation*}
For $q=0$, the same cosine sum equals $n+1$. Applying the product-to-sum formula therefore yields
\begin{align*}
&\sum_{k=1}^{n+1}
\cos\frac{\pi\ell(k-\tfrac12)}{n+1}
\cos\frac{\pi\ell'(k-\tfrac12)}{n+1}\\
&=\frac12\sum_{k=1}^{n+1}\left\{
\cos\frac{\pi(\ell-\ell')(k-\tfrac12)}{n+1}
+\cos\frac{\pi(\ell+\ell')(k-\tfrac12)}{n+1}
\right\}\\
&=\frac{n+1}{2}\mathbf1_{\{\ell=\ell'\}}.
\end{align*}
Indeed $\ell+\ell'\le2m_n<2(n+1)$, and $\ell-\ell'=0$ precisely when $\ell=\ell'$. Multiplying by the sine difference factors gives
\begin{align}
\sum_{k=1}^{n+1}
\Delta_{k}\left(\sin\frac{\pi\ell k}{n+1}\right)
\Delta_{k}\left(\sin\frac{\pi\ell' k}{n+1}\right)
&=2(n+1)\sin\frac{\pi\ell}{2(n+1)}
\sin\frac{\pi\ell'}{2(n+1)}\mathbf1_{\{\ell=\ell'\}} \nonumber \\
&=2(n+1)\sin^2\frac{\pi\ell}{2(n+1)}
\mathbf1_{\{\ell=\ell'\}} \label{eq:product_increment}.
\end{align}
Here the difference notation is extended to $k = n+1$ using its defining formula. Its final value is
\begin{align} \label{increment_sine}
\Delta_{n+1}\left(\sin\frac{\pi\ell(n+1)}{n+1}\right)
&=(-1)^\ell\sin\frac{\pi\ell}{n+1}.
\end{align}
Subtracting the product of the final differences for $\ell$ and $\ell'$ from the $(n+1)$-term sum proves the claimed $n$-term identity. The bounds follow from $|\sin x|\le |x|$. On the diagonal, the boundary term divided by $2(n+1)\sin^2(\pi\ell/(2(n+1)))$ is $2\cos^2(\pi\ell/(2(n+1)))/(n+1)$, uniformly $O(n^{-1})$. The small-angle expansion is uniform for $\ell\le m_n=o(n)$.
\end{proof}

In addition, interval integration and sine orthogonality give
\begin{align}\label{eq:gram}
\int_0^1\sin(\phi_\ell(u))\sin(\phi_{\ell'}(u))du&=\tfrac{1}{2}\mathbf1_{\{\ell=\ell'\}},\\
\int_0^1\int_0^1|\mathcal D_{n,m_n}^{j,j'}(u,s)|^2du\,ds&=m_n^{-1}.\label{eq:kernelL2}
\end{align}

\begin{proof}

Each nonzero interval has length $1/(n+1)$, so
\begin{equation*}
\int_0^1\sin(\phi_\ell(u))\sin(\phi_{\ell'}(u))du
=\frac{1}{n+1}\sum_{k=1}^{n}\sin\frac{\pi\ell k}{n+1}\sin\frac{\pi\ell'k}{n+1}.    
\end{equation*}

By the product-to-sum identity,
\begin{equation*}
    \sum_{k=1}^n
\sin\frac{\pi\ell k}{n+1}
\sin\frac{\pi\ell'k}{n+1} =
\frac12\sum_{k=1}^n\cos\frac{\pi(\ell-\ell')k}{n+1}
-\frac12\sum_{k=1}^n\cos\frac{\pi(\ell+\ell')k}{n+1}.
\end{equation*}

For $\ell\ne\ell'$, if $\ell \equiv \ell' (\mod2)$ then 
the finite geometric-sum formula therefore yields
\begin{equation*}
    \sum_{k=0}^n
\exp\left(\frac{i\pi(\ell\pm\ell')k}{n+1}\right)
=
\frac{1-\exp(i\pi(\ell\pm\ell'))}
{1-\exp(i\pi(\ell\pm\ell')/(n+1))}
= 0.
\end{equation*}

Taking real parts and subtracting the term corresponding to $k=0$, we obtain $\sum_{k=1}^n\cos\frac{\pi(\ell\pm\ell')k}{n+1}=-1$.

For $\ell \neq \ell'$, if $\ell\not\equiv\ell'\pmod 2 $ then we get
\begin{align*}
\cos\frac{\pi(\ell\pm\ell')(n+1-k)}{n+1}
&=\cos\left(\pi(\ell\pm\ell')
-\frac{\pi(\ell\pm\ell')k}{n+1}\right)\\
&=-\cos\frac{\pi(\ell\pm\ell')k}{n+1}, 
\end{align*}
and 
\begin{align*}
\sum_{k=1}^{n}\cos\frac{\pi(\ell\pm\ell')k}{n+1}
&=\sum_{k=1}^{n}\cos\frac{\pi(\ell\pm\ell')(n+1-k)}{n+1}\\
&=-\sum_{k=1}^{n}\cos\frac{\pi(\ell\pm\ell')k}{n+1}.
\end{align*}

For $\ell = \ell'$, then we obtain 
\begin{align*}
   \frac{1}{n+1} \sum_{k=1}^n
\sin\frac{\pi\ell k}{n+1}
\sin\frac{\pi\ell'k}{n+1} 
&= \frac{1}{n+1} \sum_{k=1}^n
\sin^2\frac{\pi\ell k}{n+1} \\
&= \frac{1}{2(n+1)}\left\{ \sum_{k=1}^n \cos(0) - \sum_{k=1}^n \cos\left( \frac{2\pi \ell k}{n+1} \right) \right\} \\
&= \frac{n}{2(n+1)} - \frac{-1}{2(n+1)} = \frac{1}{2}
\end{align*}
Therefore, we conclude
\begin{equation*}
    \sum_{k=1}^n
\sin\frac{\pi\ell k}{n+1}
\sin\frac{\pi\ell'k}{n+1}
=\frac{n+1}{2}\mathbf1_{\{\ell=\ell'\}},
\end{equation*}
and consequently,
\begin{equation*}
    \int_0^1\sin(\phi_\ell(u))\sin(\phi_{\ell'}(u))du
=\frac12\mathbf1_{\{\ell=\ell'\}}.
\end{equation*}

Finally expanding the finite sum in the kernel definition yields
\begin{align*}
&\int_0^1\int_0^1|\mathcal D_{n,m_n}^{j,j'}(u,s)|^2du\,ds\\
&=\frac4{m_n^2}\sum_{\ell,\ell'=1}^{m_n}
\left(\int_0^1\sin(\phi_\ell(u))\sin(\phi_{\ell'}(u))du\right)
\left(\int_0^1\sin(\phi_\ell(s))\sin(\phi_{\ell'}(s))ds\right)\\
&=\frac4{m_n^2}\sum_{\ell,\ell'=1}^{m_n}\frac14\mathbf1_{\{\ell=\ell'\}}
=\frac1{m_n}.
\end{align*}
\end{proof}

\subsection{Convergence of the signal part}
\begin{lemma}\label{Lemma_S1}
If $m_n\to\infty$ and $m_n/n\to0$, then for every $g\in C[0,1]$,

\begin{equation*}
\int_0^1\mathcal D_{n,m_n}^{j,j'}(s,s)g(s)ds\longrightarrow\int_0^1g(s)ds.
\end{equation*}
\end{lemma}

\begin{proof}
Consider the continuous kernel
\begin{align*}
    \mathcal D_{n,m_n}^{\text{cont.}}(u,s)=\frac2{m_n}\sum_{\ell=1}^{m_n}\sin(\pi\ell u)\sin(\pi\ell s).
\end{align*}

The identity $2\sin^2x=1-\cos(2x)$ and the finite cosine sum give, for $0<s<1$,
\begin{align*}
\mathcal D_{n,m_n}^{\text{cont.}}(s,s)
&=1-\frac1{m_n}\sum_{\ell=1}^{m_n}\cos(2\pi\ell s)\\
&=1+\frac1{2m_n}
-\frac{\sin((2m_n+1)\pi s)}{2m_n\sin(\pi s)}.
\end{align*}

Since
$\frac{\sin((2m_n+1)\pi s)}{\sin(\pi s)}
=1+2\sum_{\ell=1}^{m_n}\cos(2\pi\ell s),$
we have
\begin{equation*}
    \left|\frac{\sin((2m_n+1)\pi s)}{\sin(\pi s)}\right|
\le\min\left(2m_n+1,\frac1{\sin(\pi s)}\right).
\end{equation*}

Set $\delta_n=(2m_n+1)^{-1}$. The absolute quotient is symmetric about $s=1/2$, and $\sin(\pi s)\ge2s$ for $0\le s\le1/2$. Thus
\begin{align*}
\int_0^1\left|\frac{\sin((2m_n+1)\pi s)}{\sin(\pi s)}\right|ds
&\le2\left\{\int_0^{\delta_n}(2m_n+1)\,ds
+\int_{\delta_n}^{1/2}\frac1{2s}\,ds\right\}\\
&=2+\log\frac{2m_n+1}{2}
=O(\log(m_n+1)).
\end{align*}

Consequently,
\begin{align}
    \int_0^1|\mathcal D_{n,m_n}^{\text{cont.}}(s,s)-1|\,ds
&\leq \frac1{2m_n}
\left(1 + \int_0^1\left|\frac{\sin((2m_n+1)\pi s)}{\sin(\pi s)}\right|ds \right) \nonumber \\
&=O\left(\frac{\log(m_n+1)}{m_n}\right)\label{eq:continuous_diagonal_error}.
\end{align}

Next, $|\sin^2x-\sin^2y|\le|x-y|$, because the derivative of $\sin^2x$ is $\sin(2x)$. For $s\in((k-1)/(n+1),k/(n+1)]$ with $1\le k\le n$, it follows that
\begin{align*}
|\mathcal D_{n,m_n}^{j,j'}(s,s)-\mathcal D_{n,m_n}^{\text{cont.}}(s,s)|
&\le\frac2{m_n}\sum_{\ell=1}^{m_n}\pi\ell\left|\frac{k}{n+1}-s\right|\\
&\le\frac{2\pi}{m_n(n+1)}\sum_{\ell=1}^{m_n}\ell
=\frac{\pi(m_n+1)}{n+1}.
\end{align*}
On the terminal interval $(n/(n+1),1]$, the discrete diagonal is $0$ and the continuous diagonal lies between $0$ and $2$. Therefore
\begin{equation}\label{eq:discrete_diagonal_error}
\int_0^1|\mathcal D_{n,m_n}^{j,j'}(s,s)-\mathcal D_{n,m_n}^{\text{cont.}}(s,s)|\,ds
\le\frac{\pi(m_n+1)}{n+1}+\frac2{n+1}
=O\left(\frac{m_n}{n}\right).
\end{equation}
Combining \eqref{eq:continuous_diagonal_error} and \eqref{eq:discrete_diagonal_error} yields
\begin{equation*}
    \int_0^1|\mathcal D_{n,m_n}^{j,j'}(s,s)-1|\,ds
\le C\left(\frac{m_n}{n}+\frac{\log(m_n+1)}{m_n}\right),
\end{equation*}

where $C$ is independent of $n$ and $m_n$. Finally,
\begin{align*}
&\left|\int_0^1\mathcal D_{n,m_n}^{j,j'}(s,s)g(s)\,ds-\int_0^1g(s)\,ds\right|\\
&\qquad\le\|g\|_\infty\int_0^1|\mathcal D_{n,m_n}^{j,j'}(s,s)-1|\,ds\\
&\qquad\le C\|g\|_\infty
\left(\frac{m_n}{n}+\frac{\log(m_n+1)}{m_n}\right)\longrightarrow0.
\end{align*}
\end{proof}

\begin{lemma}[Signal consistency]\label{signal_consistency}
Under the assumptions above, $m_n\to\infty$ and $m_n=o(n)$ imply
$S_{n,m_n}^{j,j'}\to\int_0^1\Sigma^{j,j'}(s)ds$ in $L^2$.
\end{lemma}
\begin{proof}
It\^o's formula gives the exact expansion
\begin{align}\label{expansion}
S_{n,m_n}^{j,j'}-\int_0^1\Sigma^{j,j'}(s)ds
={}&M_{m_n}^{j,j'}(1)+M_{m_n}^{j',j}(1)\nonumber\\
&+\sum_{a=1}^3\{I_{m_n}^{a,j,j'}(1)+I_{m_n}^{a,j',j}(1)\}\nonumber\\
&+\int_0^1(\mathcal D_{n,m_n}^{j,j'}(s,s)-1)\Sigma^{j,j'}(s)ds,
\end{align}
where, for component indices $j,k$,
\begin{align*}
M_{m_n}^{j,k}(t)&=\sum_r\int_0^t\left(\sum_{r'}\int_0^s\mathcal D_{n,m_n}^{j,k}(s,u)\sigma_{r'}^k(u)dW_u^{r'}\right)\sigma_r^j(s)dW_s^r,\\
I_{m_n}^{1,j,k}(t)&=\sum_r\int_0^t\left(\int_0^s\mathcal D_{n,m_n}^{j,k}(s,u)b^k(u)du\right)\sigma_r^j(s)dW_s^r,\\
I_{m_n}^{2,j,k}(t)&=\int_0^t\left(\sum_r\int_0^s\mathcal D_{n,m_n}^{j,k}(s,u)\sigma_r^k(u)dW_u^r\right)b^j(s)ds,\\
I_{m_n}^{3,j,k}(t)&=\int_0^t\left(\int_0^s\mathcal D_{n,m_n}^{j,k}(s,u)b^k(u)du\right)b^j(s)ds.
\end{align*}
The fourth-moment assumptions, Ito isometry, BDG and Cauchy--Schwarz give
\begin{equation*}
\mathbb E|M_{m_n}^{j,k}(1)|^2+
\sum_{a=1}^3\mathbb E|I_{m_n}^{a,j,k}(1)|^2
\le C\int_0^1\int_0^1|\mathcal D_{n,m_n}^{j,k}(s,u)|^2du\,ds=O(m_n^{-1}).
\end{equation*}
BDG bounds the fourth moment of the inner stochastic integral by a constant times $(\int_0^s\mathcal D(s,u)^2du)^2$, which also handles the possibly dependent outer diffusion or drift by Cauchy--Schwarz. The diagonal term converges in $L^2$ by the deterministic bound in Lemma \ref{Lemma_S1} and $\mathbb E\sup_s|\Sigma^{j,j'}(s)|^2<\infty$.
\end{proof}

\subsection{Pure noise and cross terms}
\begin{lemma}\label{noise_bound}
Under the stated fourth noise moments,
\begin{align*}
\mathbb E[N_{n,m_n}^{j,j'}]&=\frac{(n+1)\nu^{j,j'}}{m_n}\sum_{\ell=1}^{m_n}4\sin^2\frac{\pi\ell}{2(n+1)},\\
\mathbb E|N_{n,m_n}^{j,j'}|^2&=O(m_n^4/n^2),\\
\mathbb E|N_{n,m_n}^{j,j'}-\mathbb E N_{n,m_n}^{j,j'}|^2&=O(m_n^3/n^2).
\end{align*}
If $m_n\to\infty$ and $m_n=o(n)$, the mean is
$\nu^{j,j'}\pi^2m_n^2/(3(n+1))(1+o(1))$.
\end{lemma}

\begin{proof}
Summation by parts, using the vanishing sine values at
$k=0$ and $k=n+1$, yields
\begin{equation*}
\sum_{k=1}^{n}
\sin\frac{\pi\ell k}{n+1}\,\Delta(v_k^j)
=
-\sum_{k=1}^{n+1}
\Delta_k\left(\sin\frac{\pi\ell k}{n+1}\right)v_{k-1}^j.
\end{equation*}
Consequently,
\begin{equation*}
N_{n,m_n}^{j,j'}
=\sum_{k,k'=1}^{n+1}a_{k,k'}v_{k-1}^jv_{k'-1}^{j'},
\end{equation*}
where
\begin{equation*}
a_{k,k'}
:=\frac{2}{m_n}\sum_{\ell=1}^{m_n}
\Delta_k\left(\sin\frac{\pi\ell k}{n+1}\right)
\Delta_{k'}\left(\sin\frac{\pi\ell k'}{n+1}\right).
\end{equation*}
Independence and centering imply
\begin{equation*}
\mathbb E[N_{n,m_n}^{j,j'}]
=\nu^{j,j'}\sum_{k=1}^{n+1}a_{k,k}
=\frac{4(n+1)\nu^{j,j'}}{m_n}
  \sum_{\ell=1}^{m_n}
  \sin^2\frac{\pi\ell}{2(n+1)},
\end{equation*}
where the last equality follows from
(\ref{eq:product_increment}).

To estimate the variance, use the symmetry of $a_{k,k'}$
to write
\begin{align*}
N_{n,m_n}^{j,j'}-\mathbb E N_{n,m_n}^{j,j'}
={}&
\sum_{k=1}^{n+1}a_{k,k}
\left(v_{k-1}^jv_{k-1}^{j'}-\nu^{j,j'}\right)\\
&+
\sum_{1\le k<k'\le n+1}a_{k,k'}
\left(v_{k-1}^jv_{k'-1}^{j'}
      +v_{k'-1}^jv_{k-1}^{j'}\right).
\end{align*}
Distinct centered diagonal summands are independent. For every other pair of distinct summands, a time index occurs in a single centered noise factor, so independence across time indices makes their product expectation zero. Thus all distinct summands are uncorrelated, and 
\begin{align*}
\operatorname{Var}(N_{n,m_n}^{j,j'})
={}&
\left\{
\mathbb E[(v_0^j)^2(v_0^{j'})^2]
-(\nu^{j,j'})^2
\right\}\sum_{k=1}^{n+1}a_{k,k}^2\\
&+2\left\{
\nu^{j,j}\nu^{j',j'}
+(\nu^{j,j'})^2
\right\}
\sum_{1\le k<k'\le n+1}a_{k,k'}^2.
\end{align*}
Here the coefficient in the first line is the variance of a diagonal noise product. The coefficient in the second line follows from
\begin{equation*}
    \mathbb E\left[
\left(
v_{k-1}^jv_{k'-1}^{j'}
+v_{k'-1}^jv_{k-1}^{j'}
\right)^2
\right]
=
2\left\{
\nu^{j,j}\nu^{j',j'}
+(\nu^{j,j'})^2
\right\},
\qquad k<k'.
\end{equation*}

Finite fourth moments make both coefficients finite. Since the noise distribution does not depend on $n$, and $a_{k,k'}=a_{k',k}$, we obtain
\begin{equation*}
    \operatorname{Var}(N_{n,m_n}^{j,j'})
\le C\sum_{k,k'=1}^{n+1}a_{k,k'}^2,
\end{equation*}
where $C$ is independent of $n$ and $m_n$.

By \eqref{eq:product_increment},
\begin{align*}
\sum_{k,k'=1}^{n+1}a_{k,k'}^2
&=\frac{4}{m_n^2}
\sum_{\ell,\ell'=1}^{m_n}
\left\{
\sum_{k=1}^{n+1}
\Delta_k\left(\sin\frac{\pi\ell k}{n+1}\right)
\Delta_k\left(\sin\frac{\pi\ell'k}{n+1}\right)
\right\}^2\\
&=\frac{16(n+1)^2}{m_n^2}
\sum_{\ell=1}^{m_n}
\sin^4\frac{\pi\ell}{2(n+1)}\\
&\le
\frac{\pi^4}{m_n^2(n+1)^2}
\sum_{\ell=1}^{m_n}\ell^4
=O\left(\frac{m_n^3}{n^2}\right).
\end{align*}
This proves the centered second-moment bound.
Similarly, the mean formula and $|\sin x|\le |x|$ give
$|\mathbb E N_{n,m_n}^{j,j'}|=O(m_n^2/n)$. Hence
\begin{equation*}
\mathbb E[|N_{n,m_n}^{j,j'}|^2]
=
\operatorname{Var}(N_{n,m_n}^{j,j'})
+|\mathbb E N_{n,m_n}^{j,j'}|^2
=O\left(\frac{m_n^4}{n^2}\right).
\end{equation*}

Finally, if $m_n=o(n)$, then uniformly for $\ell\le m_n$,
\begin{equation*}
4\sin^2\frac{\pi\ell}{2(n+1)}
=
\frac{\pi^2\ell^2}{(n+1)^2}
\left\{1+O\left(\frac{m_n^2}{(n+1)^2}\right)\right\}.
\end{equation*}
Substituting this into the mean formula and using $\sum_{\ell=1}^{m_n}\ell^2
\sim\frac{m_n^3}{3}$
proves the final assertion. 
\end{proof}

\begin{lemma}\label{cross_bound}
The cross term satisfies 
\begin{equation*}
\mathbb E|C_{n,m_n}^{j,j'}|^2=O(m_n/n).
\end{equation*}
\end{lemma}

\begin{proof}
Put
\begin{equation*}
\xi_\ell^j=\sum_{k=1}^n r_{k,\ell}^n\Delta X_k^j,
\qquad
\eta_\ell^j=\sum_{k=1}^n r_{k,\ell}^n\Delta(v_k^j).
\end{equation*}
Cauchy--Schwarz, Ito isometry and
$\sum_{k=1}^n(r_{k,\ell}^n)^2=1$ yield
\begin{equation*}
\mathbb E[(\xi_\ell^j)^2]\le\frac{C}{n+1}.
\end{equation*}

Summation by parts gives
\begin{equation*}
\eta_\ell^j
=-\sqrt{\frac{2}{n+1}}
\sum_{k=1}^{n+1}
\Delta_k\left(\sin\frac{\pi\ell k}{n+1}\right)v_{k-1}^j.
\end{equation*}
Here the term $k=n+1$ incorporates the original terminal
term, since
\begin{equation*}
\Delta_{n+1}\left(
\sin\frac{\pi\ell(n+1)}{n+1}\right)
=-\sin\frac{\pi\ell n}{n+1}.
\end{equation*}

By independence, centering of the noise vectors and \eqref{eq:product_increment},
\begin{equation*}
\mathbb E[\eta_\ell^j\eta_{\ell'}^{j'}]
=
4\nu^{j,j'}\sin^2\frac{\pi\ell}{2(n+1)}
\mathbf1_{\{\ell=\ell'\}}.
\end{equation*}

Since the noise is independent of the signal,
\begin{align*}
\mathbb E\left[
\left|\frac{n+1}{m_n}
\sum_{\ell=1}^{m_n}\xi_\ell^j\eta_\ell^{j'}\right|^2
\right]
&=\frac{4(n+1)^2\nu^{j',j'}}{m_n^2}
\sum_{\ell=1}^{m_n}
\mathbb E[(\xi_\ell^j)^2]
\sin^2\frac{\pi\ell}{2(n+1)}\\
&\le\frac{C(n+1)}{m_n^2}
\sum_{\ell=1}^{m_n}
\sin^2\frac{\pi\ell}{2(n+1)}\\
&\le\frac{C}{m_n^2(n+1)}
\sum_{\ell=1}^{m_n}\ell^2
=O\left(\frac{m_n}{n}\right).
\end{align*}
Finally,
\begin{equation*}
C_{n,m_n}^{j,j'}
=\frac{n+1}{m_n}
\sum_{\ell=1}^{m_n}
\left(\xi_\ell^j\eta_\ell^{j'}
+\xi_\ell^{j'}\eta_\ell^j\right).
\end{equation*}
Applying the same estimate with $j$ and $j'$ interchanged,
and using $(a+b)^2\le2a^2+2b^2$, proves the assertion.
\end{proof}

\begin{theorem}[Consistency]\label{thm:consistency}
Under the assumptions in Section 2.1, if $m_n\to\infty$ and $m_n=o(n^{1/2})$, then
\begin{equation*}
\hat V_{n,m_n}^{j,j'}\xrightarrow{L^2}\int_0^1\Sigma^{j,j'}(s)ds.
\end{equation*}
\end{theorem}
\begin{proof}
Combine the exact decomposition with Lemmas \ref{signal_consistency}, \ref{noise_bound} and \ref{cross_bound}.
\end{proof}

\section{Asymptotic normality}\label{sec:normality}
\subsection{First-order Malliavin assumptions and signal limit}
Retain the model, filtration, and moment assumptions of Section \ref{sec:consistency}. For asymptotic normality, add the following differentiability assumption; standard notation for the Malliavin derivative and its Sobolev spaces is as in \cite{ND}.
\begin{assume}[First-order Malliavin regularity]\label{ass:malliavin}
For every finite $p\ge2$, each $\sigma_r^j$ belongs to $\mathbb D^{1,p}(L^2[0,1])$. Its time sections belong to $\mathbb D^{1,p}$ for almost every time. There are jointly measurable derivatives $D_u^q\sigma_r^j(t)$ such that
\begin{equation*}
    \mathbb E\left[\operatorname*{ess\,sup}_{u,t\in[0,1]}
|D_u^q\sigma_r^j(t)|^p\right]<\infty,
\qquad 1\le q,r\le d,\quad 1\le j\le J.
\end{equation*}
For fixed $u$, use adapted versions of the derivative processes in $t$, vanishing for $u>t$. Essential suprema are over the deterministic time variables.
\end{assume}
No Malliavin differentiability of $b$ is imposed.
Fix $j,j'$ and write $m=m_n$ within the proof. Set
\begin{equation*}
    \mathcal M_n^{j,j'}(t):=M_m^{j,j'}(t)+M_m^{j',j}(t),\qquad
\mathcal V^{j,j'}(t):=\int_0^t
\{\Sigma^{j,j}(s)\Sigma^{j',j'}(s)+(\Sigma^{j,j'}(s))^2\}ds.
\end{equation*}

The processes $M_m^{j,k}$ and $I_m^{a,j,k}$ are those in \eqref{expansion}.
\begin{theorem}[Signal stable limit under Malliavin regularity]\label{prop:signal_clt}
Under the assumptions of Section \ref{sec:consistency} and Assumption \ref{ass:malliavin}, if $m_n\to\infty$ and $m_n=o(n^{2/3})$, then
\begin{equation*}
    \sqrt{m_n}\left(S_{n,m_n}^{j,j'}-\int_0^1\Sigma^{j,j'}(s)ds\right)
\xrightarrow{st.}\mathcal{MN}(0,\mathcal V^{j,j'}(1))
\end{equation*}
relative to $\mathcal F_1$. The limit is conditionally centered Gaussian with the displayed variance, on an extension of the probability space. 
\end{theorem}
We prove the theorem in the following subsections.

\subsection{Deterministic kernel bounds and concentration}
For brevity, write $m=m_n$ and
$f_\ell(u)=\sin(\phi_\ell(u))$, and set
\begin{equation*}
K_n(s,u):=\mathcal D_{n,m_n}^{j,j'}(s,u)
=\frac2m\sum_{\ell=1}^m f_\ell(s)f_\ell(u).
\end{equation*}
For $u\in(t_{k-1},t_k]$, $1\le k\le n$, write
$q_n(u)=t_k$.

\begin{lemma}[Kernel bounds]\label{mall:lem:kernel}
For $1\le m<n$,
\begin{align}
\int_0^1K_n(s,u)^2du &= \frac{K_n(s,s)}{m} \leq \frac{2}{m},\label{mall:eq:row}
\end{align}
More generally, let
\begin{equation*}
    \mathcal{K}_n^{(p)}=\max\left\{\sup_s\int_0^1|K_n(s,u)|^pdu,
\sup_u\int_0^1|K_n(s,u)|^pds\right\}, 
\end{equation*}
then 
\begin{align}
    \mathcal{K}_n^{(p)} \leq \begin{cases}
        C_1\log(m+1)/m \qquad p=1, \\
        C_p/m \qquad p>1.
    \end{cases}
\end{align}
Also
\begin{equation}\label{mall:eq:diag}
\int_0^1|K_n(s,s)-1|ds
\le C\left(\frac mn+\frac{\log(m+1)}m\right).
\end{equation}
\end{lemma}
For short, put $\mathcal{K}_n = \mathcal{K}_n^{(1)}$. 
\begin{proof}
Expanding the finite sum and using orthogonality, we obtain
\begin{align*}
\int_0^1 K_n(s,u)^2du
&=\frac4{m^2}\sum_{\ell,\ell'=1}^{m}
f_\ell(s)f_{\ell'}(s)
\int_0^1f_\ell(u)f_{\ell'}(u)du\\
&=\frac2{m^2}\sum_{\ell=1}^{m}f_\ell(s)^2
=\frac{K_n(s,s)}m
\le\frac2m.
\end{align*}

Define $C_m(x):=\frac1m\sum_{\ell=1}^m\cos(\pi\ell x)$.
The product-to-sum
identity gives
\begin{equation*}
K_n(s,u)
=C_m(q_n(s)-q_n(u))-C_m(q_n(s)+q_n(u)),
\end{equation*}
for $s,u$ in the observed intervals.
The finite geometric-sum formula yields
\begin{equation*}
|C_m(x)|
\le
\min\left\{
1,\frac1{m\,\operatorname{dist}(x,2\mathbb Z)}
\right\},
\end{equation*}
with the bound interpreted as $1$ at zero distance.
Indeed, for $x\notin2\mathbb Z$,
\begin{equation*}
\left|\sum_{\ell=1}^m e^{i\pi\ell x}\right|
\le\frac1{|\sin(\pi x/2)|}
\le\frac1{\operatorname{dist}(x,2\mathbb Z)}.
\end{equation*}

For $a,b\in[0,1]$, 
\begin{equation*}
\operatorname{dist}(a-b,2\mathbb Z)=|a-b|,
\quad
\operatorname{dist}(a+b,2\mathbb Z)
=\min\{a+b,2-a-b\}\ge|a-b|.
\end{equation*}
We consequently get
\begin{equation*}
|K_n(s,u)|
\le\frac2{m|q_n(s)-q_n(u)|}
\end{equation*}
whenever $q_n(s)\ne q_n(u)$.
Also, $|K_n(s,u)|\le2$ everywhere.

On the observed intervals,
$0\le q_n(s)-s\le(n+1)^{-1}$, so
\begin{equation*}
|q_n(s)-q_n(u)|
\ge |s-u|-\frac1{n+1}.
\end{equation*}
Since $m<n$, whenever $|s-u|\ge2/m$ we obtain
\begin{equation*}
|q_n(s)-q_n(u)|\ge\frac{|s-u|}{2},
\qquad
|K_n(s,u)|\le\frac4{m|s-u|}.
\end{equation*}
The latter kernel bound also holds when either variable
is outside the observed intervals, since the kernel is
then zero.

For $m\ge3$, set
$I_s:=\{u\in[0,1]:|s-u|<2/m\}$.
Splitting the integral over $I_s$ and its complement gives
\begin{equation*}
\int_0^1|K_n(s,u)|du
\le \int_{I_s} 2 du +
\int_{I_s^c}
 \frac4m \frac{du}{|s-u|}
\le \frac{C\log(m+1)}m.
\end{equation*}
For $m=1,2$, the same conclusion follows from
$|K_n|\le2$, after increasing $C$.
Symmetry gives the identical bound for the column
integrals. Hence
\begin{equation*}
\mathcal K_n\le \frac{ C\log(m+1)}m,
\end{equation*}
where $C$ is independent of $n$ and $m$.
For $p>1$, the same decomposition, together with
$|K_n(s,u)|\le2$ and the preceding off-diagonal bound,
gives
\begin{equation*}
\int_0^1|K_n(s,u)|^pdu
\le \frac{2^{p+2}}m
+\frac{2\cdot4^p}{m^p}\int_{2/m}^1r^{-p}dr
\le\frac{C_p}{m},\qquad p>1, 
\end{equation*}
uniformly in $s\in[0,1]$ and $m \geq 3$. For $m=1,2$, the same bound follows from $|K_n|\le2$
after increasing $C_p$. Symmetry of the kernel gives the corresponding column bound.

Finally, the diagonal estimate was proved in Section 2
by comparison with the continuous sine kernel.
Since $K_n=\mathcal D_{n,m_n}^{j,j'}$, that estimate gives
\begin{equation*}
\int_0^1|K_n(s,s)-1|ds
\le C\left(\frac mn+\frac{\log(m+1)}m\right).
\end{equation*}
\end{proof}

\begin{lemma}[Weighted concentration]\label{mall:lem:concentration}
If $m\to\infty$ and $m=o(n)$, then for every continuous $g$ on $[0,1]^2$ and every $t\in[0,1]$,
\begin{equation}\label{mall:eq:concentration}
m\int_0^t\int_0^sK_n(s,u)^2g(s,u)duds
\longrightarrow\frac12\int_0^tg(s,s)ds.
\end{equation}
The same assertion holds pathwise for any random continuous $g$. For every $\delta>0$,
\begin{equation}\label{mall:eq:offdiag}
\sup_s m\int_0^{(s-\delta)_+}K_n(s,u)^2du\longrightarrow0.
\end{equation}
\end{lemma}

\begin{proof}
By the proof of
Lemma \ref{mall:lem:kernel}, for every fixed $\delta>0$
and all sufficiently large $n$,
\begin{equation*}
\sup_{|s-u|\ge\delta}|K_n(s,u)|
\le\frac{C}{m\delta}.
\end{equation*}
Consequently,
\begin{equation*}
\sup_{s\in[0,1]}
m\int_{\substack{u\in[0,1]\\|s-u|\ge\delta}}
K_n(s,u)^2\,du
\le\frac{C}{m\delta^2}\longrightarrow0.
\end{equation*}
This proves \eqref{mall:eq:offdiag}.
Fix $g\in C([0,1]^2)$ and $t\in[0,1]$, and set
\begin{equation*}
\bar g(s,u):=\frac{g(s,s)+g(u,u)}2.
\end{equation*}
Let $\omega_g$ denote the modulus of continuity of $g$
with respect to the maximum norm.
When $|s-u|<\delta$,
$|g(s,u)-\bar g(s,u)|\le\omega_g(\delta)$.
Using the preceding off-diagonal estimate and
$m\iint_{[0,1]^2}K_n(s,u)^2\,du\,ds=1$, we obtain
\begin{equation*}
\left|
m\int_0^t\int_0^s
K_n(s,u)^2\{g(s,u)-\bar g(s,u)\}\,du\,ds
\right|
\le
\omega_g(\delta)
+\frac{C\|g\|_\infty}{m\delta^2}.
\end{equation*}
Taking $n\to\infty$ and then $\delta\downarrow0$
shows that this difference vanishes.

By symmetry,
\begin{align*}
m\int_0^t\int_0^sK_n(s,u)^2\bar g(s,u)\,du\,ds
&=\frac m2\int_0^t\int_0^t
K_n(s,u)^2g(s,s)\,du\,ds\\
&=\frac12\int_0^tK_n(s,s)g(s,s)\,ds\\
&\quad-\frac m2\int_0^t\int_t^1
K_n(s,u)^2g(s,s)\,du\,ds,
\end{align*}
where the last equality uses \eqref{mall:eq:row}.
Splitting the outer integral of the boundary term
at $(t-\delta)_+$ gives
\begin{equation*}
m\int_0^t\int_t^1K_n(s,u)^2\,du\,ds
\le\frac{C}{m\delta^2}
+\int_{(t-\delta)_+}^tK_n(s,s)\,ds
\le\frac{C}{m\delta^2}+2\delta.
\end{equation*}
Thus the boundary term vanishes by taking $n\to\infty$
and then $\delta\downarrow0$.
Moreover, \eqref{mall:eq:diag} implies
\begin{equation*}
\left|
\int_0^t\{K_n(s,s)-1\}g(s,s)\,ds
\right|
\le
\|g\|_\infty\int_0^1|K_n(s,s)-1|\,ds
\longrightarrow0.
\end{equation*}
Combining these estimates proves
\eqref{mall:eq:concentration}.

For a random weight with almost surely continuous sample
paths, apply the deterministic argument to each sample
path on the event of continuity.
\end{proof}

\subsection{Negligibility of the drift remainders}
\begin{lemma}[Drift remainders]\label{mall:lem:drift}
If $m\to\infty$ and $m=o(n)$, then $\sqrt m I_m^{a,j,k}(1)\to0$ in $L^2$, and hence in probability, for $a=1,2,3$.
\end{lemma}
\begin{proof}
Recall $\mathcal{K}_n \le C\frac{\log(m+1)}m$.
It\^o isometry and the coefficient moment assumptions give
\begin{align*}
\mathbb E|I_m^{1,j,k}(1)|^2
&\le\mathbb E\left[\sup_{s\le1}|b^k(s)|^2\sup_{s \leq 1}\sum_r|\sigma_r^j(s)|^2\right]
\int_0^1 \left( \int_0^s|K_n(s,u)|du \right)^2ds\le C\mathcal{K}_n^2,\\
\mathbb E|I_m^{3,j,k}(1)|^2
&\le\mathbb E[\sup_{s\le1}|b^j(s)|^2\sup_{s\le1}|b^k(s)|^2]
\left(\int_0^1  \int_0^s|K_n(s,u)|du ds\right)^2\le C\mathcal{K}_n^2.
\end{align*}
Consequently $m\mathbb E|I_m^{a,j,k}(1)|^2\le C\log^2(m+1)/m\to0$ for $a=1,3$.

For $a=2$ abbreviate
\begin{equation*}
Y_\sigma(v,s)=\sum_r\int_0^vK_n(s,u)\sigma_r^k(u)dW_u^r,
\qquad v\le s.
\end{equation*}
BDG and \eqref{mall:eq:row} imply, for every finite $p\ge2$,
\begin{equation*}
\|Y_\sigma(v,s)\|_p\le C_p\left(\int_0^vK_n(s,u)^2du\right)^{1/2}
\le C_p m^{-1/2}.
\end{equation*}
Since $I_m^{2,j,k}(1)=\int_0^1Y_{\sigma}(s,s)b^j(s)ds$, Fubini and symmetry give
\begin{equation*}
m\mathbb E|I_m^{2,j,k}(1)|^2
=2m\int_0^1\int_0^s
\mathbb E[Y_\sigma(s,s)b^j(s)Y_\sigma(t,t)b^j(t)]dtds.
\end{equation*}
Fix $\delta\in(0,1)$. On the strip $(s-\delta)_+\le t\le s$, H\"older with four $L^4$ factors bounds the absolute integrand by $C/m$. The strip contribution, including the prefactor $2m$, is therefore at most $C\delta$.

On $t<s-\delta$, necessarily $s>\delta$, split
\begin{equation*}
Y_\sigma(s,s)=Y_\sigma(s-\delta,s)+R_{s,\delta},\qquad
R_{s,\delta}:=\sum_r\int_{s-\delta}^sK_n(s,u)\sigma_r^k(u)dW_u^r.
\end{equation*}
For the old part, the same H\"older estimate and BDG give the absolute contribution at most
\begin{equation*}
C\int_\delta^1\left(m\int_0^{s-\delta}K_n(s,u)^2du\right)^{1/2}ds\longrightarrow0
\end{equation*}
by \eqref{mall:eq:offdiag}. For the recent part,
$\mathbb E[R_{s,\delta}\mid\mathcal F_{s-\delta}]=0$, and
$b^j(s-\delta)Y_\sigma(t,t)b^j(t)$ is $\mathcal F_{s-\delta}$-measurable. Hence
\begin{equation*}
\mathbb E[R_{s,\delta}b^j(s)Y_\sigma(t,t)b^j(t)] =\mathbb E[R_{s,\delta}\{b^j(s)-b^j(s-\delta)\}
Y_\sigma(t,t)b^j(t)].
\end{equation*}
Since $\|R_{s,\delta}\|_4+\|Y_\sigma(t,t)\|_4\le Cm^{-1/2}$, H\"older yields
\begin{equation*}
m\left|\mathbb E[R_{s,\delta}\{b^j(s)-b^j(s-\delta)\}
Y_\sigma(t,t)b^j(t)]\right|
\le C\|b^j(s)-b^j(s-\delta)\|_4.
\end{equation*}
Integration over $t <s-\delta$ and then H\"older in $s$ give
\begin{equation*}
\limsup_n m\mathbb E|I_m^{2,j,k}(1)|^2
\le C\delta+C\left(\mathbb E\int_\delta^1
|b^j(s)-b^j(s-\delta)|^4ds\right)^{1/4}.
\end{equation*}
Extend $b^j$ by zero outside $[0,1]$. Its fourth-moment bound implies membership in $L^4(\Omega\times\mathbb R)$, whose time translations are strongly continuous. Thus the last term tends to zero as $\delta\downarrow0$. Taking $n\to\infty$ first and then $\delta\downarrow0$ proves $m\mathbb E|I_m^{2,j,k}(1)|^2\to0$. 
\end{proof}
By \eqref{mall:eq:diag}, the last term of \eqref{expansion}, after multiplication by $\sqrt m$, is
\begin{equation}\label{mall:eq:diagonal_clt}
O_p\left(m^{3/2}/n+\log(m+1)/\sqrt m\right)=o_p(1)
\end{equation}
if $m\to\infty$ and $m=o(n^{2/3})$.

\subsection{Moment bounds and Malliavin duality estimates}
We first derive moment bounds and Malliavin derivative formulas for the inner stochastic integrals associated with general coefficients. These estimates are then used to establish two convergence lemmas needed for the signal martingale. In particular,
\begin{equation}\label{mall:eq:smallquantities}
m \mathcal{K}_n^2+\sqrt m\,\mathcal{K}_n\le C\{\log^2(m+1)/m+\log(m+1)/\sqrt m\}\longrightarrow0.
\end{equation}
Let $\vartheta=(\vartheta_1,\ldots,\vartheta_d)$ be any progressively measurable vector coefficient with the coefficient moment bounds of Section \ref{sec:consistency} and the first-derivative bounds of Assumption \ref{ass:malliavin}, and let $c(s)$ be an adapted scalar coefficient with the same bounds. We will apply the estimates to $\vartheta=\sigma^j$, $\vartheta=\sigma^j+\sigma^{j'}$, and to $c=\sigma_r^j$ or $c=\Sigma^{j,k}$. Products defining the latter satisfy the required bounds by the first-order chain rule and H\"older.
Define, for deterministic $s$ and $v\le s$,
\begin{align*}
Y_{\vartheta}(v,s)&=\sum_r\int_0^vK_n(s,w)\vartheta_r(w)dW_w^r,\\
V_{\vartheta}(v,s)&=\sum_r\int_0^vY_{\vartheta}(w,s)K_n(s,w)\vartheta_r(w)dW_w^r.
\end{align*}
For every finite $p\ge2$, BDG and H\"older give the uniform estimates
\begin{equation}\label{mall:eq:YVbounds}
\|Y_{\vartheta}(v,s)\|_p\le C_p\sqrt{(2m_n^{-1})},\qquad
\|V_{\vartheta}(v,s)\|_p\le C_p (2m_n^{-1}).
\end{equation}
For almost every $u$ and every Wiener index $q$ the first-derivative formulas are
\begin{align}
&D_u^qY_{\vartheta}(v,s)
=\mathbf1_{\{u\le v\}}K_n(s,u)\vartheta_q(u)
+\sum_r\int_0^vK_n(s,w)D_u^q \vartheta_r(w)dW_w^r,\label{mall:eq:DY}\\
&D_u^qV_{\vartheta}(v,s)
=\mathbf1_{\{u\le v\}}Y_{\vartheta}(u,s)K_n(s,u)\vartheta_q(u)
+\sum_r\int_0^vK_n(s,w)D_u^q(Y_{\vartheta}(w,s)\vartheta_r(w))dW_w^r.
\label{mall:eq:DV}
\end{align}
The derivatives of adapted coefficients vanish before $u$. Thus the displayed derivative integrands, for fixed $u$, are adapted. These formulas follow by closure from smooth adapted approximations in the stated Malliavin Sobolev spaces.
In particular, define
\begin{equation*}
J_{s,u}^q:=\sum_r\int_0^sK_n(s,v)D_u^q(Y_{\vartheta}(v,s)\vartheta_r(v))dW_v^r.
\end{equation*}
The product rule and \eqref{mall:eq:DY} give the useful bound
\begin{equation}\label{mall:eq:Jbound}
\|J_{s,u}^q\|_p
\le C_p\left(\mathbf1_{\{u\le s\}}|K_n(s,u)|\sqrt{(2m_n^{-1})}+(2m_n^{-1})\right).
\end{equation}
The direct term in \eqref{mall:eq:DY} contributes
$C_p\mathbf1_{\{u\le s\}}|K_n(s,u)|m^{-1/2}$.
The remaining terms are bounded by $C_pm^{-1}$
using the product rule, H\"older, BDG, and
\eqref{mall:eq:row}.

\begin{lemma}[Asymptotic orthogonality]\label{mall:lem:orthog}
For each fixed $t$,
\begin{equation*}
\mathbb E\left|\sqrt m\int_0^tY_{\vartheta}(s,s)c(s)ds\right|^2
\le C(m \mathcal{K}_n^2+\sqrt m\,\mathcal{K}_n)\longrightarrow0.
\end{equation*}
The resulting processes converge to zero uniformly in probability on $[0,1]$.
\end{lemma}
\begin{proof}
Write $G_n(s,s')=\mathbb E[Y_{\vartheta}(s,s)Y_{\vartheta}(s',s')c(s)c(s')]$. Fubini gives 
\begin{equation*}
    \mathbb E|\sqrt m\int_0^tY_{\vartheta}(s,s)c(s)ds|^2=m\int_{[0,t]^2}G_n(s,s')dsds'.
\end{equation*}
By Malliavin duality,
\begin{equation*}
G_n(s,s')=\sum_q\int_0^sK_n(s,u)
\mathbb E\left[\vartheta_q(u)D_u^q\{Y_{\vartheta}(s',s')c(s)c(s')\}\right]du.
\end{equation*}
The product rule and \eqref{mall:eq:DY} split $G_n(s,s')$ into $E_1+E_2+E_3$, where
\begin{align*}
E_1={}&\sum_q\int_0^{s\wedge s'}K_n(s,u)K_n(s',u)
\mathbb E[\vartheta_q(u)^2c(s)c(s')]du,\\
E_2={}&\sum_q\int_0^sK_n(s,u)
\mathbb E\left[\vartheta_q(u)c(s)c(s')L_{s',u}^q\right]du,\\
L_{s',u}^q:={}&\sum_r\int_0^{s'}K_n(s',v)D_u^q\vartheta_r(v)dW_v^r,\\
E_3={}&\sum_q\int_0^sK_n(s,u)
\mathbb E\left[\vartheta_q(u)Y_{\vartheta}(s',s')D_u^q(c(s)c(s'))\right]du.
\end{align*}
The coefficient moment assumptions imply
$|E_1|\le C\int_0^{s\wedge s'}|K_n(s,u)K_n(s',u)|du$.
Changing the order of the deterministic integrals gives
\begin{align*}
\int_{[0,t]^2}|E_1|dsds'
&\le C\int_0^t\left(\int_u^t|K_n(s,u)|ds\right)^2du
\le C\mathcal{K}_n^2.
\end{align*}
For $E_2$, BDG and the derivative moment assumptions give, for each finite $p\ge2$,
\begin{equation*}
\|L_{s',u}^q\|_p\le C_p\left(\int_0^{s'}K_n(s',v)^2dv\right)^{1/2}
\le C_p\sqrt{(2m_n^{-1})}.
\end{equation*}
H\"older therefore yields
$|E_2|\le C\sqrt{(2m_n^{-1})}\int_0^s|K_n(s,u)|du$,
and integration gives $\int_{[0,t]^2}|E_2|dsds'\le C\sqrt{(2m_n^{-1})}\mathcal{K}_n$.
For $E_3$, use
\begin{equation*}
D_u^q(c(s)c(s'))=(D_u^qc(s))c(s')+c(s)D_u^qc(s')
\end{equation*}
together with $\|Y_{\vartheta}(s',s')\|_p\le C_p\sqrt{(2m_n^{-1})}$ and H\"older. This gives the same bound
$\int_{[0,t]^2}|E_3|dsds'\le C\sqrt{(2m_n^{-1})}\mathcal{K}_n$.
Consequently
\begin{equation*}
\mathbb E|\sqrt m\int_0^tY_{\vartheta}(s,s)c(s)ds|^2\le Cm(\mathcal{K}_n^2+\sqrt{(2m_n^{-1})}\mathcal{K}_n)
\le C(m\mathcal{K}_n^2+\sqrt m\,\mathcal{K}_n)\longrightarrow0.
\end{equation*}
The first estimate uses the column $L^1$ bound; the other two use the row $L^1$ bound. 
For uniform convergence, Cauchy--Schwarz in time and H\"older give, for $0\le v\le t\le1$,
\begin{align*}
\mathbb E\left|\sqrt m \left\{\int_0^tY_{\vartheta}(s,s)c(s)ds-\int_0^vY_{\vartheta}(s,s)c(s)ds \right\}\right|^2
&\le m(t-v)\int_v^t\mathbb E|Y_{\vartheta}(s,s)c(s)|^2ds\\
&\le Cm(2m_n^{-1})(t-v)^2\le C(t-v)^2.
\end{align*}
The continuous processes start at zero and are tight in $C[0,1]$ by the Kolmogorov tightness criterion. Their fixed-time limits are zero; hence they converge to the zero process in distribution in $C[0,1]$. Since the limit is deterministic, this is uniform convergence in probability.
\end{proof}

\begin{lemma}[Quadratic-variation remainder]\label{mall:lem:bracketerror}
For each fixed $t$,
\begin{equation}\label{mall:eq:bracketerror}
\mathbb E\left|m\int_0^t V_{\vartheta}(s,s)c(s)ds\right|^2
\le C(m \mathcal{K}_n^2+\sqrt m\,\mathcal{K}_n)\longrightarrow0.
\end{equation}
\end{lemma}
\begin{proof}
 A single use of duality on the integral defining $V_{\vartheta}(s,s)$ gives
\begin{equation*}
\mathbb E[V_{\vartheta}(s,s)V_{\vartheta}(s',s')c(s)c(s')] =\sum_q\int_0^s K_n(s,u)
\mathbb E\left[Y_{\vartheta}(u,s)\vartheta_q(u)
D_u^q\{V_{\vartheta}(s',s')c(s)c(s')\}\right]du.
\end{equation*}
Decompose the right-hand side into $F_1,F_2,F_3$, where $F_1$ uses the direct term in \eqref{mall:eq:DV}, $F_2$ uses $J_{s',u}^q$, and $F_3$ differentiates $c(s)c(s')$. Explicitly,
\begin{align*}
F_1={}&\sum_q\int_0^{s\wedge s'}K_n(s,u)K_n(s',u) \mathbb E[Y_{\vartheta}(u,s)Y_{\vartheta}(u,s')\vartheta_q(u)^2c(s)c(s')]du,\\
F_2={}&\sum_q\int_0^sK_n(s,u)
\mathbb E[Y_{\vartheta}(u,s)\vartheta_q(u)J_{s',u}^q c(s)c(s')]du,\\
F_3={}&\sum_q\int_0^sK_n(s,u)
\mathbb E[Y_{\vartheta}(u,s)\vartheta_q(u)V_{\vartheta}(s',s')D_u^q(c(s)c(s'))]du.
\end{align*}
By \eqref{mall:eq:YVbounds} and H\"older,
\begin{equation*}
\int_{[0,t]^2}|F_1|dsds'\le C (2m_n^{-1}) \mathcal{K}_n^2,
\qquad
\int_{[0,t]^2}|F_3|dsds'\le C (2m_n^{-1})^{3/2}\mathcal{K}_n.
\end{equation*}
For $F_2$, use \eqref{mall:eq:Jbound} and $\|Y_{\vartheta}(u,s)\|_p\le C_p\sqrt{(2m_n^{-1})}$, obtaining
\begin{equation*}
|F_2|\le C\int_0^s|K_n(s,u)|
\left\{(2m_n^{-1})|K_n(s',u)|\mathbf1_{\{u\le s'\}}+(2m_n^{-1})^{3/2}\right\}du.
\end{equation*}
Its double integral is bounded by $C((2m_n^{-1}) \mathcal{K}_n^2+(2m_n^{-1})^{3/2}\mathcal{K}_n)$. Multiplying these estimates by $m^2$, and using $(2m_n^{-1})=2/m$, proves \eqref{mall:eq:bracketerror}. This is an $L^2$ estimate for the random remainder, not merely convergence of its expectation.
\end{proof}

\subsection{Quadratic variation and asymptotic orthogonality}
Write
\begin{equation*}
Z_n^j(s)=\sum_r\int_0^sK_n(s,u)\sigma_r^j(u)dW_u^r.
\end{equation*}
Then we obtain
\begin{align*}
\mathcal M_n^{j,j'}(t)
=\sum_r\int_0^t\{Z_n^{j'}(s)\sigma_r^j(s)
+Z_n^j(s)\sigma_r^{j'}(s)\}dW_s^r,
\end{align*}
and
\begin{align}\label{mall:eq:bracket}
m\langle\mathcal M_n^{j,j'}\rangle_t
=m\int_0^t\{&(Z_n^{j'}(s))^2\Sigma^{j,j}(s)
+2Z_n^{j'}(s)Z_n^j(s)\Sigma^{j,j'}(s)\nonumber\\
&+(Z_n^j(s))^2\Sigma^{j',j'}(s)\}ds.
\end{align}
For the following weighted-concentration step, take $\vartheta$ and $c$ to have continuous sample paths, in addition to the preceding moment and first-derivative bounds. This continuity is satisfied by every coefficient used below, since $\vartheta$ is a diffusion vector or a sum of diffusion vectors and $c$ is a covariance entry. If $m\to\infty$ and $m=o(n)$, then for each fixed $t$,
\begin{equation}\label{mall:eq:weighted-square}
m\int_0^tY_{\vartheta}(s,s)^2c(s)ds
\xrightarrow{p}\frac12\int_0^t|\vartheta(s)|^2c(s)ds.
\end{equation}
\begin{proof}[Proof of \eqref{mall:eq:weighted-square}]
Hold $s$ fixed. As a process in $v\le s$, $Y_{\vartheta}(v,s)$ is a continuous martingale starting at zero, with
\begin{equation*}
\langle Y_{\vartheta}\rangle_v=\int_0^vK_n(s,u)^2|\vartheta(u)|^2du.
\end{equation*}
It\^o's formula applied to $Y_{\vartheta}(v,s)^2$, and then evaluated at $v=s$, gives
\begin{align*}
Y_{\vartheta}(s,s)^2
&=\int_0^sK_n(s,u)^2|\vartheta(u)|^2du
+2\sum_r\int_0^sY_{\vartheta}(u,s)K_n(s,u)\vartheta_r(u)dW_u^r\\
&=\int_0^sK_n(s,u)^2|\vartheta(u)|^2du+2V_{\vartheta}(s,s).
\end{align*}
The moment bounds justify integration of this identity in $s$. Thus
\begin{align*}
m\int_0^tY_{\vartheta}(s,s)^2c(s)ds&=T_{1,n}(t)+T_{2,n}(t),\\
T_{1,n}(t)&:=m\int_0^t\int_0^sK_n(s,u)^2|\vartheta(u)|^2c(s)duds,\\
T_{2,n}(t)&:=2m\int_0^tV_{\vartheta}(s,s)c(s)ds.
\end{align*}
On the probability-one event where $\vartheta$ and $c$ are continuous, the function
$g_\omega(s,u)=|\vartheta(u,\omega)|^2c(s,\omega)$ is a continuous function on $[0,1]^2$. Lemma \ref{mall:lem:concentration} applies to each such sample path, yielding
\begin{equation*}
T_{1,n}(t)\longrightarrow\frac12\int_0^t|\vartheta(s)|^2c(s)ds
\quad\text{almost surely}.
\end{equation*}
This is an application to a random continuous weight by fixing the sample path; no expectation or independence is used here.
For the second term, Lemma \ref{mall:lem:bracketerror} gives
\begin{equation*}
\mathbb E|T_{2,n}(t)|^2
\le C\bigl(m\mathcal{K}_n^2+\sqrt m\,\mathcal{K}_n\bigr)
\le C\left\{\frac{\log^2(m+1)}m+
\frac{\log(m+1)}{\sqrt m}\right\}\longrightarrow0.
\end{equation*}
Consequently $T_{2,n}(t)\to0$ in probability. Combining the almost sure convergence of $T_{1,n}(t)$ with this $L^2$ convergence proves \eqref{mall:eq:weighted-square}.

\end{proof}
Choosing $\vartheta=\sigma^{j'}$, $c=\Sigma^{j,j}$ gives half of the product-variance integral. The choice $\vartheta=\sigma^j$, $c=\Sigma^{j',j'}$ gives the other half. Polarization of \eqref{mall:eq:weighted-square}, using $\vartheta=\sigma^j+\sigma^{j'}$ and $c=\Sigma^{j,j'}$, gives
\begin{equation*}
m\int_0^tZ_n^j(s)Z_n^{j'}(s)\Sigma^{j,j'}(s)ds
\xrightarrow{p}\frac12\int_0^t(\Sigma^{j,j'}(s))^2ds.
\end{equation*}
Consequently \eqref{mall:eq:bracket} converges at each fixed time in probability to $\mathcal V^{j,j'}(t)$. The left-hand side is an increasing process and the limit is continuous and increasing. A finite-grid monotonicity argument therefore upgrades convergence to uniform convergence in probability on $[0,1]$.

For each Wiener component,
\begin{equation*}
\sqrt m\langle\mathcal M_n^{j,j'},W^r\rangle_t
=\sqrt m\int_0^t\{Z_n^{j'}(s)\sigma_r^j(s)
+Z_n^j(s)\sigma_r^{j'}(s)\}ds.
\end{equation*}
Lemma \ref{mall:lem:orthog} proves uniform convergence to zero in probability. To obtain the condition for every bounded signal martingale $H$, use its predictable Wiener representation
$H=H_0+\sum_r\int h_r\,dW^r$.
Approximate $h$ in $L^2(\Omega\times[0,1])$ by bounded elementary predictable processes. For each such approximation, the covariation is a finite sum of bounded coefficients times increments of the preceding Wiener covariations and hence converges uniformly to zero. The approximation error is bounded pathwise by
\begin{equation*}
\sup_{t\le1}|\langle\sqrt m\,\mathcal M_n^{j,j'},
\textstyle\sum_r\int(h_r-h_r^{(L)})dW^r\rangle_t|
\le \{m\langle\mathcal M_n^{j,j'}\rangle_1\}^{1/2}
\left(\int_0^1|h-h^{(L)}|^2ds\right)^{1/2}.
\end{equation*}
The first factor is tight by the proved bracket convergence; the second tends to zero in probability as $L\to\infty$. Thus the required orthogonality holds for all bounded signal martingales.

\subsection{Completion of the proof and the noisy-data stable limit}

The continuous-martingale stable limit theorem \cite{JJ} gives
\begin{equation*}
\sqrt m\,\mathcal M_n^{j,j'}(1)
\xrightarrow{st.}
\mathcal{MN}(0,\mathcal V^{j,j'}(1))
\end{equation*}
relative to $\mathcal F_1$.
Lemma \ref{mall:lem:drift}, \eqref{mall:eq:diagonal_clt},
and \eqref{expansion}, together with stable Slutsky,
therefore complete the proof of Theorem \ref{prop:signal_clt}.
In particular, the negligibility of the normalized drift
and diagonal remainders, bracket convergence, and asymptotic
orthogonality have been derived from the assumptions of
Section \ref{sec:consistency}, Assumption \ref{ass:malliavin},
and the deterministic kernel estimates.

We now pass from the signal estimator to the noisy estimator.
Set
\begin{equation*}
\mathcal G
:=\mathcal F_1\vee\sigma(v_k:k\ge0).
\end{equation*}
 By independence of the signal and noise, the signal stable limit extends from $\mathcal F_1$ to $\mathcal G=\mathcal F_1^v$.
This follows first for products of bounded signal- and noise-measurable test variables and then for all bounded $\mathcal G$-measurable test variables by $L^1$ approximation.

\begin{theorem}[Asymptotic normality under Malliavin regularity]
\label{thm:normality_noise}
Under the assumptions of Section \ref{sec:consistency}
and Assumption \ref{ass:malliavin}, if $m_n\to\infty$
and $m_n=o(n^{2/5})$, then
\begin{equation*}
\sqrt{m_n}\left(
\hat V_{n,m_n}^{j,j'}
-\int_0^1\Sigma^{j,j'}(s)\,ds
\right)
\xrightarrow{st.}
\mathcal{MN}(0,\mathcal V^{j,j'}(1))
\end{equation*}
relative to $\mathcal G$.
\end{theorem}

\begin{proof}
The exact decomposition gives
\begin{align*}
\sqrt{m_n}\left(
\hat V_{n,m_n}^{j,j'}
-\int_0^1\Sigma^{j,j'}(s)\,ds
\right)
&=\sqrt{m_n}
\left(S_{n,m_n}^{j,j'}-\int_0^1\Sigma^{j,j'}(s)\,ds\right) \\
& \qquad +\sqrt{m_n}C_{n,m_n}^{j,j'}
+\sqrt{m_n}N_{n,m_n}^{j,j'}.
\end{align*}
By Lemmas \ref{noise_bound} and \ref{cross_bound},
\begin{equation*}
\mathbb E\left|
\sqrt{m_n}N_{n,m_n}^{j,j'}
\right|^2
=O(m_n^5/n^2),
\qquad
\mathbb E\left|
\sqrt{m_n}C_{n,m_n}^{j,j'}
\right|^2
=O(m_n^2/n).
\end{equation*}
Both bounds tend to zero when $m_n=o(n^{2/5})$.
Moreover, this bandwidth condition implies
$m_n=o(n^{2/3})$, as required by
Theorem \ref{prop:signal_clt}.
The signal stable limit relative to $\mathcal G$,
established above, and stable Slutsky complete the proof.
\end{proof}

\section{Feasible bias correction}\label{sec:bias}
\subsection{Estimating the noise covariance}
All assumptions and symbols are those of Section \ref{sec:consistency}. In particular, initial noise is included and may be correlated across components. Define
\begin{equation}\label{eq:nuhat}
\hat\nu_n^{j,j'}=\frac1{2n}\sum_{k=1}^n\Delta Y_k^j\Delta Y_k^{j'},
\qquad \hat\nu_n=(\hat\nu_n^{j,j'})_{j,j'=1}^J.
\end{equation}
\begin{prop}[Noise covariance estimation]\label{prop:nuhat}
For each fixed pair $j,j'$, $\mathbb E|\hat\nu_n^{j,j'}-\nu^{j,j'}|^2\le C/n$. Consequently $\hat\nu_n^{j,j'}-\nu^{j,j'}=O_p(n^{-1/2})$. The matrix $\hat\nu_n$ is symmetric positive semidefinite.
\end{prop}
\begin{proof}
 The coefficient moment assumptions, H\"older and BDG imply, uniformly in $k$,
\begin{equation*}
\mathbb E|\Delta X_k^j|^2\le C(n+1)^{-1},\qquad
\mathbb E|\Delta X_k^j|^4\le C(n+1)^{-2}.
\end{equation*}
Indeed the fourth moments of the drift and diffusion increments are bounded by $C(n+1)^{-4}$ and $C(n+1)^{-2}$, respectively. Expand
\begin{align*}
\hat\nu_n^{j,j'}-\nu^{j,j'}={}&\frac1{2n}\sum_{k=1}^n
\{\Delta(v_k^j)\Delta(v_k^{j'})-2\nu^{j,j'}\}\\
&+\frac1{2n}\sum_{k=1}^n\Delta X_k^j\Delta(v_k^{j'})
+\frac1{2n}\sum_{k=1}^n\Delta(v_k^j)\Delta X_k^{j'}\\
&+\frac1{2n}\sum_{k=1}^n\Delta X_k^j\Delta X_k^{j'}.
\end{align*}
Let $\mathcal H_n:=\sigma(\Delta X_1^j,\ldots,\Delta X_n^j)$.
By independence of the signal and noise,
\begin{equation*}
\mathbb E[v_k^{j'}v_h^{j'}\mid\mathcal H_n]
=\nu^{j',j'}\mathbf1_{\{k=h\}},
\qquad 0\le k,h\le n.
\end{equation*}
Moreover, summation by parts gives
\begin{equation*}
\sum_{k=1}^n\Delta X_k^j\Delta(v_k^{j'})
=
-\Delta X_1^jv_0^{j'}
+\sum_{k=1}^{n-1}
(\Delta X_k^j-\Delta X_{k+1}^j)v_k^{j'}
+\Delta X_n^jv_n^{j'}.
\end{equation*}
Since the signal increments are $\mathcal H_n$-measurable
and the noise is independent of $\mathcal H_n$,
taking conditional second moments in the preceding
summation-by-parts identity gives
\begin{align*}
\mathbb E\left[
\left|\sum_{k=1}^n\Delta X_k^j\Delta(v_k^{j'})\right|^2
\Bigm|\mathcal H_n
\right]
& =
\nu^{j',j'}\left\{
|\Delta X_1^j|^2+|\Delta X_n^j|^2
+\sum_{k=1}^{n-1}
|\Delta X_k^j-\Delta X_{k+1}^j|^2
\right\}\\
&\quad\le4\nu^{j',j'}\sum_{k=1}^n|\Delta X_k^j|^2.
\end{align*}
Taking expectations and using
$\mathbb E|\Delta X_k^j|^2\le C/(n+1)$ yields
\begin{equation*}
\mathbb E\left|
\frac1{2n}\sum_{k=1}^n
\Delta X_k^j\Delta(v_k^{j'})
\right|^2
\le \frac{4\nu^{j',j'}}{4n^2}\sum_{k=1}^n
\mathbb E|\Delta X_k^j|^2
= O\left( n^{-2} \right).
\end{equation*}
The other mixed term satisfies the same bound
with $j$ and $j'$ interchanged.

Put $U_k:=\Delta(v_k^j)\Delta(v_k^{j'})-2\nu^{j,j'}$. The centering and temporal independence of the noise vectors imply
$\mathbb E U_k=0$.
Moreover, $U_k$ depends only on $v_{k-1}$ and $v_k$,
so $(U_k)$ is $1$-dependent.
The finite fourth noise moments give
$\sup_k\mathbb E U_k^2<\infty$.
Consequently, only diagonal and adjacent terms contribute
to the second moment of the sum, and Cauchy--Schwarz yields
\begin{equation*}
\mathbb E\left|\frac1{2n}\sum_{k=1}^nU_k\right|^2 =\frac1{4n^2}
\left\{
\sum_{k=1}^n\mathbb E U_k^2
+2\sum_{k=1}^{n-1}\mathbb E[U_kU_{k+1}]
\right\} \le\frac{C}{n},
\end{equation*}
where $C$ is independent of $n$. Finally, 
\begin{equation*}
\left\|
\frac1{2n}\sum_{k=1}^n
\Delta X_k^j\Delta X_k^{j'}
\right\|_2\le
\frac1{2n}\sum_{k=1}^n
\|\Delta X_k^j\|_4
\|\Delta X_k^{j'}\|_4 \le \frac C{n+1}.
\end{equation*}
The signal contribution produces only a small mean error:
\begin{equation*}
\mathbb E\hat\nu_n^{j,j'}-\nu^{j,j'}
=\frac1{2n}\sum_{k=1}^n\mathbb E[\Delta X_k^j\Delta X_k^{j'}]=O(n^{-1}).
\end{equation*}
Combining the four terms proves the bound; Chebyshev gives the probabilistic order. Positive semidefiniteness follows from the outer-product representation in \eqref{eq:nuhat}.
\end{proof}

\subsection{Oracle and feasible corrections}
Define the deterministic coefficient
\begin{equation}\label{eq:beta}
\beta_{n,m_n}:=\frac{4(n+1)}{m_n}\sum_{\ell=1}^{m_n}
\sin^2\frac{\pi\ell}{2(n+1)}.
\end{equation}
Lemma \ref{noise_bound} gives $\mathbb EN_{n,m_n}^{j,j'}=\beta_{n,m_n}\nu^{j,j'}$. Also
$\beta_{n,m_n}\le C m_n^2/n$ for $m_n<n$, and
$\beta_{n,m_n}\sim\pi^2m_n^2/(3(n+1))$ if $m_n\to\infty$, $m_n=o(n)$.
The oracle correction (known noise covariance) and the feasible correction are
\begin{align}
\hat V_{n,m_n}^{j,j',\mathrm{or}}
&:=\hat V_{n,m_n}^{j,j'}-\beta_{n,m_n}\nu^{j,j'},\label{eq:oracle}\\
\hat V_{n,m_n}^{j,j',\mathrm{bc}}
&:=\hat V_{n,m_n}^{j,j'}-\beta_{n,m_n}\hat\nu_n^{j,j'}.
\label{bc}
\end{align}
Their exact decompositions, using the same signal notation as Sections 2--3, are
\begin{align}
\hat V_{n,m_n}^{j,j',\mathrm{or}}
&=S_{n,m_n}^{j,j'}+C_{n,m_n}^{j,j'}+
(N_{n,m_n}^{j,j'}-\mathbb EN_{n,m_n}^{j,j'}),\nonumber\\
\hat V_{n,m_n}^{j,j',\mathrm{bc}}
&=\hat V_{n,m_n}^{j,j',\mathrm{or}}
-\beta_{n,m_n}(\hat\nu_n^{j,j'}-\nu^{j,j'}).
\label{eq:bc_decomposition}
\end{align}


\subsection{Limit theorems for the corrected estimators}
\begin{theorem}[Corrected consistency]\label{thm:bc_consistency}
Under the assumptions of Section \ref{sec:consistency}, if $m_n\to\infty$ and $m_n=o(n^{2/3})$, both corrected estimators converge in $L^2$, and hence in probability, to $\int_0^1\Sigma^{j,j'}(s)ds$.
\end{theorem}
\begin{proof}
Lemma \ref{signal_consistency} applies because $m_n=o(n)$. Using \eqref{eq:bc_decomposition}, Lemmas \ref{noise_bound}--\ref{cross_bound}, and Proposition \ref{prop:nuhat}, we obtain
\begin{align*}
\mathbb E\left|\hat V_{n,m_n}^{j,j',\mathrm{bc}}-
\int_0^1\Sigma^{j,j'}(s)ds\right|^2
\le C\bigg\{&\mathbb E\left|S_{n,m_n}^{j,j'}-
\int_0^1\Sigma^{j,j'}(s)ds\right|^2\\
&+\frac{m_n}{n}+\frac{m_n^3}{n^2}+\frac{m_n^4}{n^3}\bigg\}.
\end{align*}
All terms vanish under the stated condition. For the oracle estimator the last term is absent.
\end{proof}
 We next establish a stable central limit theorem
under a smaller bandwidth range.

\begin{theorem}[Corrected stable limit]\label{thm:bc_clt}
Under the assumptions of Section \ref{sec:consistency} and Assumption \ref{ass:malliavin}, let $m_n\to\infty$ and $m_n=o(n^{1/2})$. Then
\begin{align*}
&\sqrt{m_n}\left(\hat V_{n,m_n}^{j,j',\text{or}}-\int_0^1\Sigma^{j,j'}(s)ds\right)
\xrightarrow{st.}\mathcal{MN}(0,\mathcal V^{j,j'}(1)), \\
&\sqrt{m_n}\left(\hat V_{n,m_n}^{j,j',\text{bc}}-\int_0^1\Sigma^{j,j'}(s)ds\right)
\xrightarrow{st.}\mathcal{MN}(0,\mathcal V^{j,j'}(1)),
\end{align*}
with the same meaning of stable convergence and mixed normality as in Section \ref{sec:normality}.
\end{theorem}
\begin{proof}
By Lemmas \ref{noise_bound}--\ref{cross_bound},
\begin{equation*}
\mathbb E|\sqrt{m_n}(N_{n,m_n}^{j,j'}-\mathbb EN_{n,m_n}^{j,j'})|^2
=O(m_n^4/n^2),\qquad
\mathbb E|\sqrt{m_n}C_{n,m_n}^{j,j'}|^2=O(m_n^2/n).
\end{equation*}
Both tend to zero. Moreover, Proposition \ref{prop:nuhat} and \eqref{eq:beta} give
\begin{equation*}
\mathbb E\left|\sqrt{m_n}\beta_{n,m_n}
(\hat\nu_n^{j,j'}-\nu^{j,j'})\right|^2
\le C\frac{m_n^5}{n^3}\longrightarrow0.
\end{equation*}
Theorem \ref{prop:signal_clt}, \eqref{eq:bc_decomposition}, and stable Slutsky prove both limits.
\end{proof}

\section{Conclusion}
For the sine-based estimator introduced in Section 1 under common equidistant sampling, we establish consistency when $m_n\to\infty$ and $m_n=o(n^{1/2})$. The pure-noise mean is $O(m_n^2/n)$, its centered second moment is $O(m_n^3/n^2)$, and the second moment of the full cross term is $O(m_n/n)$. Consequently, $m_n=o(n^{2/5})$ is sufficient for both noise contributions to vanish at the $\sqrt{m_n}$ scale.

For adapted diffusion coefficients satisfying the first-order Malliavin assumptions, stable asymptotic normality follows from the verified drift remainders, kernel concentration, bracket convergence, and asymptotic orthogonality. No derivative of the drift or second Malliavin derivative of the diffusion coefficient is required. With estimated noise covariance, the feasible bias correction is consistent for $m_n\to\infty$, $m_n=o(n^{2/3})$, and has a centered stable mixed-normal limit for $m_n=o(n^{1/2})$ under the same Malliavin assumptions. These sufficient bandwidth conditions do not establish a general optimal-rate stable limit at $m_n\asymp\sqrt n$.

\end{document}